\documentclass[12pt,oneside,draft]{amsart}
\usepackage{amsfonts,amssymb,amscd,amsthm,amsmath}
\usepackage{graphicx}
\usepackage[english]{babel}


\newtheorem{theorem}[subsection]{Theorem}

\newtheorem{proposition}[subsection]{Proposition}
\newtheorem{lemma}[subsection]{Lemma}

\theoremstyle{definition}
\newtheorem{definition}[subsection]{Definition}
\theoremstyle{remark}

\newcommand{\N}{{\mathcal N}}
\newcommand{\mt}[1]{\operatorname{#1}}

\newcommand{\OO}{{\mathcal O}}

\newcommand{\ZZ}{{\mathbb Z}}
\newcommand{\CC}{{\mathbb C}}

\newcommand{\PP}{{\mathbb P}}

\newcommand{\NN}{{\mathbb N}}

\newcommand{\cyc}[1]{\ZZ_{#1}}

\title{Divisorial contractions to some Gorenstein singularities}
\author{I. Yu. Fedorov}

\date{}
\email{ifedorov@mi.ras.ru}
\begin{document}
\begin{abstract} Divisorial contractions to singularities, defined by equations
$xy+z^n+u^n=0$ $n\ge3$ and $xy+z^3+u^4=0$ are classified.
\end{abstract}
\maketitle
 The problem of  birational classification of algebraic
varieties is highly interconnected with the problem of description
of singularities on them. One of the most important class of
three-dimensional singularities is terminal singularities, which
arise within minimal models programm. Despite the analytical
classification of the singularities
\cite{Dan82},\cite{Reid83},\cite{Ms84},\cite{Mori85}, this
description does not help one to fully understand many birationl
properties of them. In particular, the problem of description of
resolution of such singularities and the problem of classification
of morphisms of terminal varieties are still up-to-date.
Divisorial contractions to cyclic quotient singularities were
described by Y.Kawamata \cite{Kaw1}, S.Mori \cite{Mori} and S.
Cutkosky \cite{Cut} classified contractions from terminal
Gorenstein threefolds. T.Luo \cite{Luo} set out contractions when
the index is not increase. A.Corti \cite{Corti} with M.Mella \cite{CorMe}
described divisorial contractions to $xy+z^n+u^n=0$ points, where $n=2,3$.
Recently M.Kawakita \cite{Kawk0}, \cite{Kawk1}, \cite{Kawk2} has gave a
description of contractions to a smooth and $cA$ points.
In this paper we classify divisorial contractions from a terminal 3-folds
to a germ of a point defined by the equation $xy+z^n+u^n=0$, where  $n\ge3$ and to a germ of a singularity
defined by the equation $xy+z^3+u^4=0$ using quite different method then the one
introduced in \cite{Kawk2}. Our method
allow us to deal with all terminal Gorenstein singularities an with non Gorenstein of a type $cA/m$.
\par
The author would like to thank Professor V.A. Iskovskikh and
Professor Yu.G. Prokhorov for their vulnerable  discussions and
encouragement. The author was partially supported by grants RFBR-99-01-01132, RFBR-96-15-96146  and
INTAS-OPEN-97-2072.

\section{\bf Preliminary results}
We will deal with varieties over $\CC$. The basic results and
notions are contained in \cite{KMM}, \cite{Sh1}.
\begin{definition}
Consider a cyclic quotient singularity $X:=\CC^n/\cyc{m}(a_1,
\dots,a_n)$, where $ a_i\in\NN$ and $\gcd(a_1,\dots,a_n)=1$ (the
case $m=1$, i.~e. $X\simeq\CC^n $, is also possible). Let $x_1,
\dots, x_n$ be eigen-coordinates in $\CC^n $, for $\cyc{m}$.
\textit{The weighted blow-up}\index{weighted blow-up} of $X$ with
weights $a_1,\dots,a_n$ is a projective birational morphism
$f\colon Y\to X$ such that $Y$ is covered by affine charts
$U_1,\dots,U_n$, where
\[
\begin{array}{ccc}
U_i=\CC^n_{y_1,\dots,y_n}/\cyc{a_i}(-a_1,\dots,&m,&\dots,-a_n).\\
&\uparrow&\\&i&\\
\end{array}
\]
The coordinates in $X$ and in $U_i$ are related by
\[
x_i=y_i^{a_i/m}, \qquad x_j=y_jy_i^{a_j/m},\quad j\ne i.
\]
The exceptional set $E$ of $f$ is an irreducible divisor and
$E\cap U_i=\{y_i=0 \}/\cyc{a_i}$. The morphism $f\colon Y\to X$ is
toric, i.e. there is an equivariant natural action of $(\CC^*)^n$.
It is easy to show that $E$ is the weighted projective space
$\PP(a_1,\dots, a_n)$ and $\OO_E(bE)=\OO_{\PP}(-mb)$, if $b$ is
divisible by $\mt{lcm}(a_1,\cdots,a_n)$ (and then $bE$ is a Cartier
divisor).

Note that the blow-up constructed above depends on a choice of
numbers $a_1,\dots,a_n$, and not just on their values $\mod m$.

Let $X$ be a hypersurface in $\CC^n$. By {\it weighted blow-up} of $X$
with weights $(a_1,\dots,a_n)$ we will mean the restriction  of the weighted blow-up
of $\CC^n$ with weights $(a_1,\dots,a_n)$ on $X$.
\end{definition}
\section{\bf Contractions to $xy+z^n+u^n=0$}
In this section we will prove the theorem.
\begin{theorem}
Let $f:Y\to X$ be a divisorial contractionfrom from terminal 3-fold $Y$ to $X$ -- a germ of $xy+z^n+u^n=0$, $n\ge 3$
singularity, $\rho(Y/X)=1$, divisor $S$ is $f$-exceptional. Then $f$ is isomorphic
to the weighted blow-up of $X$ with weights $(k,n-k,1,1)$ for some integer $1\le k\le n-1$.
\end{theorem}

The classification we will obtain using the following plan:
Let $f:Y\to X$ be a divisorial contractionfrom from terminal 3-fold $Y$ to a germ of $xy+z^n+u^n=0$, $n\ge 3$
singularity $X$, such that $f$-exceptional divisor $S$ is an irreducible reduced divisor.
Let $f_1:Y_1\to X$ be the weighted blow-up with weights $(1,n-1,1,1)$. It follows from \cite{Fed2}
that the discrepansy of the $f_1$-exceptional divisor $E_1$ is equal to $1$.
\par
In fact, $E_1$ is a rational surface with just one singularity of  $\frac 1{n-1}(1,1)$ type, $\rho(S)$
is equal to $n$, the configuration of  $(-1)$-curves $\bar {l_i}$, $i=1,\dots,n$ on the minimal resolution
of $\bar E_1\to E_1$ is as on the Fig. 1.
\par
We will deal with following cases:
\begin{enumerate}
\item The discrepancy of $a(S)$ (over $X$) is equal to 1. All those morphisms (in Mori's cathegory)
is classified in  \cite{Fed2}.

\item $a(S)=k\ge 2$.
\begin{enumerate}
\item The center of $S$ on $E_1$ is a point and it does not lie on any $l_i$.
\item The center of $S$ on $E_1$ is a curve and it does not coincide with any $l_i$.
\item The center of $S$ on $E_1$ is a point and it lies on some $l_i$.
\item The center of $S$ on $E_1$ is a curve and it coincides with some $l_i$.
\end{enumerate}
\end{enumerate}
We will prove our main theorem checking all the cases.

\subsection{Geometry of $E_1$}
\begin{proposition}
For the surface $E_1$ the following statements are true:
\begin{enumerate}
\item $E_1$ is a rational surface with one singularity of  $\frac
1{n-1}(1,1)$ type;
\item the Picard number of $E_1$ is equal to $n$;

\item the configuration of $(-1)$-curves $\bar {l_i}$, $i=1,\dots,n$ on the minimal resolution of $E_1$
$\alpha:\bar{E_1}\to E_1$ is as on the Fig. 1.
\end{enumerate}
\end{proposition}
\begin{proof}
\begin{enumerate}
\item We have $Y_1$ covered by four affine charts $U_1,\ U_2,\\ U_3,\ U_4$:
\begin{gather*}
U_1=\{\bar y+\bar z^{n}+\bar u^{n}=0\}\subseteq\CC^4,\\
U_2=\{\bar x+\bar z^{n}+\bar u^{n}=0\}/\ZZ_{n-1}(-1,1,-1,-1),\\
U_3=\{\bar x\bar y+1+\bar u^{n}=0\}\subseteq\CC^4,\\
U_4=\{\bar x\bar y+\bar z^{n}+1=0\}\subseteq\CC^4.
\end{gather*}
Thus, there is only one singularity of $\frac 1 {n-1}(-1,1,1)$ type on $Y_1$. It lies in chart $U_2$.
\item We have  $K_{E_1}^2=\frac{(n-n+2)^2n}{n-1}=\frac{4n}{n-1}$. The minimal resolution of $E_1$ is just one
blow-up, and the exceptionl curve will be  $-(n-1)$-curve with discrepancy $\frac{3-n}{n-1}$. Hence
\begin{gather*}
K_{\bar E_1}^2=K_{E_1}^2+\left(\frac{3-n}{n-1}\right)(-2+n-1)=
\frac{4n}{n-1}-\\-\frac{9+n^2-6n}{n-1} =\frac{10n-9-n^2}{n-1}=9-n
\end{gather*}

It follows from the Noether formulae that $\chi(\bar E_1)=n+3$. Thus, we have $\rho(E_1)=n$.

\item Let's look at the divisor $(x=0)$ on $E_1$. It consists of $n$ curves $l_1,\dots,l_n$.

Let's prove that on the $\bar E_1$ $(\bar l_i)^2=-1$:

We have $(x=0)|_{E_1}=l_1+\dots+l_n$. Self-intersection numbers of $l_i$ are equal since the symmetry of $l_i$.
We have
\begin{equation*}
\frac n{n-1}=(x=0)^2|_{E_1}=(\sum l_i)^2=nl_i^2+\sum_{i<j} 2l_il_j=nl_i^2+n
\end{equation*}

Therefore
\begin{equation*}
l_i^2=\frac 1{n-1}-1=-\left(\frac{n-2}{n-1}\right)
\end{equation*}
and
\begin{equation*}
\bar l_i^2=l_i^2-\frac 1{n-1}=-1.
\end{equation*}
\end{enumerate}
\end{proof}
\subsection{Weighted blow-ups of $xy+z^n+u^n=0$}
\begin{lemma}
There are no weighted blow-ups of $xy+z^n+u^n=0$ producing an irreducible divisor with discrepancy $k\ge 2$
in Mori's category.
\end{lemma}
\proof
Let's look at weighted blow-up of $h:Y\to X$ with weights $(a,b,c,d)$. With no loss of generality
we will consider that $a+b=nc$. We will consider two cases:
\begin{enumerate}
\item $a+b=nc=nd$
\item $a+b=nc<nd$
\end{enumerate}
\begin{enumerate}
\item Let's examine the terminality of $Y$. Indeed, $Y$ is covered  by four affine charts $U_1,\ U_2,\ U_3,\ U_4$
\begin{gather*}
U_1=\{\bar y+\bar z^n+\bar u^n=0\}/\ZZ_a(1,-b,-c,-c),\\
U_2=\{\bar x+\bar z^n+\bar u^n=0\}/\ZZ_b(-a,1,-c,-c),\\
U_3=\{\bar x\bar y+\bar u^{n}+1=0\}\ZZ_c(-a,-b,1,0),\\
U_4=\{\bar x\bar y+\bar z^{n}+1=0\}\ZZ_c(-a,-b,0,1).
\end{gather*}
Classification of terminal singularities (see \cite{YPG}) tells us that singularities
of $Y$ in charts $U_1$ and $U_2$ are terminal if and only if
\begin{gather*}
c\equiv 1\mod a\\
c\equiv 1\mod b.
\end{gather*}
Then, either $c=1$ or $c>a$ and $c>b$. The latter is impossible since $a+b=nc$.
\item  In this case we have $d>c$. Let's check the terminality of $Y$ again.
It is covered by four affine charts  $U_1,\ U_2,\ U_3,\ U_4$
\begin{gather*}
U_1=\{\bar y+\bar z^n+\bar x^{n(d-c)}\bar u^n=0\}/\ZZ_a(1,-b,-c,-d),\\
U_2=\{\bar x+\bar z^n+\bar y^{n(d-c)}\bar u^n=0\}/\ZZ_b(-a,1,-c,-d),\\
U_3=\{\bar x\bar y+\bar z^{n(d-c)}\bar u^{n}+1=0\}\ZZ_c(-a,-b,1,-d),\\
U_4=\{\bar x\bar y+\bar z^{n}+\bar u^{n(d-c)}=0\}\ZZ_d(-a,-b,-c,1).
\end{gather*}
Then, singularities in chart $U_4$ are terminal if $d>c$. \qed
\end{enumerate}
\par Here are some results from the paper \cite{Fed2}.
\begin{theorem}
Let $X$ be a germ of a 3-dimentional terminal $cA$ point defined by $xy+f(z,u)=0$,
$n$ ia a number of divisors with discrepancy 1 over $X$.
Then, we have that $n=deg_{min}(f)-1$, where $deg_{min}(f)$ is minimal degree among
degrees of all the monoms in $f$.
\end{theorem}
Moreover, it was showed in that paper that if we make the weighted blow-up of $X$ with weights $(a,b,1,1)$
then the others divisors with discrepancy one lies precisely over non Gorenstein points of this weighted blow-up.
\subsection{Examination of the cases}
\subsubsection{Cases 2(a) and 2(b)}

Looking at the different models (among all the weighted blow-ups with weights $(a,b,1,1)$)
we can reach the situation when the center of $S$ lies in the singularity of the model. We obviously
can consider only $Y_1$ and $Y_{n-1}$. It follows from \cite{Kaw1}, that  $S$ can be realized
by some weighted blow-up $\bar Y\to X$. It is easy to notice that in this case $\bar Y \simeq Y$.
Lemma 3.2 produces a contradiction with a terminality of $Y$.

\subsubsection{Cases 2(c) and 2(d)}

It follows from the \cite{Fed2}, that all terminal varieties which realize all the divisors with discrepancy 1 over $X$
differs from one to other in flops in $\bar l_i$.
The exact consequence
\begin{equation*}
0\longrightarrow\OO_{P^1}(-1)\longrightarrow\N_{X|l_i}
\longrightarrow\OO_{P^1}(-1)\longrightarrow 0
\end{equation*}
allow us to blow-up  $\bar  l_i$ and then contract the surface to another ruling (see \cite{Kul}).
Therefore, we can consider the center of $S$ on some model to be in the  $x^2+y^2+z^2+u^2=0$.
It follows from \cite{Corti} that in this case $S$ realizes by the ordinary blow-up of this point.
Hence, in our model $S$ realizes as a blow-up of $\bar l_i$. Thus, the case 2(c) is not possible.
\par

The case 2(d) we will consider on the weighted blow-up of $X$ with weights $(1,n-1,1,1)$,
$E_1$ is the exceptional divisor.
It follows from \cite{Kaw1} that a divisorial contraction to a curve passing through
a terminal cyclic quotient singularity is a weighted blow-up of this singularity. Therefore, we get a
contradiction with lemma 2.5.
\qed

The main theorem of this section is proved.

\section{\bf Contractions to $xy+z^3+u^4=0$}
In this section we will prove the theorem.
\begin{theorem}
Let $f:Y\to X$ be a divisorial contractionfrom from terminal 3-fold $Y$ to $X$ -- a germ of $xy+z^3+u^4=0$
singularity, $\rho(Y/X)=1$, divisor $S$ is $f$-exceptional. Then $f$ is isomorphic
to the weighted blow-up of $X$ with weights $(k,n-k,1,1)$ for $k=1,2$.
\end{theorem}

We will use the same method and the same plan that used in the previous section.

\begin{enumerate}
\item The discrepancy of $a(S)$ (over $X$) is equal to 1. All those morphisms (in Mori's cathegory)
is classified in  \cite{Fed2}.

\item $a(S)=k\ge 2$.
\begin{enumerate}
\item The center of $S$ on $E_1$ is a point and it does not lie on any $l_i$.
\item The center of $S$ on $E_1$ is a curve and it does not coincide with any $l_i$.
\item The center of $S$ on $E_1$ is a point and it lies on some $l_i$.
\item The center of $S$ on $E_1$ is a curve and it coincides with some $l_i$.
\end{enumerate}
\end{enumerate}
We will prove our main theorem checking all the cases.

\subsection{Geometry of $E_1$}
\begin{proposition}
For the surface $E_1$ the following statements are true:
\begin{enumerate}
\item $E_1$ is a rational surface with two Du Val singularities of a type $A_1, A_2$;
\item the Picard number of $E_1$ is equal to $1$;
\item there is just one $(-1)$-curve $\bar {l}$ on the minimal resolution of $E_1$
$\alpha:\bar{E_1}\to E_1$ and $l$ is passing through two singularities of $E_1$.
\end{enumerate}
\end{proposition}
\begin{proof}
\begin{enumerate}
\item We have $Y_1$ covered by four affine charts $U_1,\ U_2,\\ U_3,\ U_4$:
\begin{gather*}
U_1=\{\bar y+\bar z^{3}+\bar u^{4}=0\}\subseteq\CC^4,\\
U_2=\{\bar x+\bar z^{3}+\bar u^{4}=0\}/\ZZ_{2}(-1,1,-1,-1),\\
U_3=\{\bar x\bar y+1+\bar u^{4}=0\}\subseteq\CC^4,\\
U_4=\{\bar x\bar y+\bar z^{3}+\bar u=0\}\subseteq\CC^4.
\end{gather*}
There is only one singularity on $Y_1$ of a $\frac 12(-1,1,1)$ type. It lies in the chart $U_2$.
In the chart $U_4$ there is a singularity of a type $A_2$ on $E_1$.

\item We have $K_{E_1}^2=\frac{(3-3+2)^23}{2}=6$. The minimal resolution of $E_1$ is consist of
3 blow-ups. We have,$K_{\bar E_1}^2=K_{E_1}^2=6$

It follows from the Noether formulae that  $\chi(\bar E_1)=4$. Therefore, $\rho(E_1)=1$.

\item Actually, there is only one surface with such a properties up to an isomorphism.
This is $\PP(1,2,3)$. There is only one $l$ such that on the minimal resolution $\bar l^2=-1$ on it.
\end{enumerate}
\end{proof}
\subsection{Weighted blow-ups of $xy+z^3+u^4=0$}
\begin{lemma}
There are no weighted blow-ups of $xy+z^3+u^4=0$ producing an
irreducible divisor with discrepancy $k\ge 2$ in Mori's category.
\end{lemma}
\proof
Let's look at weighted blow-up of $h:Y\to X$ with weights $(a,b,c,d)$. We will consider the following three cases:
\begin{enumerate}
\item $a+b=3c=4d$;
\item $a+b=3c<4d$;
\item $a+b=4d<3c$.
\end{enumerate}
\begin{enumerate}
\item Let's examine the terminality of $Y$. Indeed, $Y$ is covered  by four affine charts $U_1,\ U_2,\ U_3,\ U_4$
\begin{gather*}
U_1=\{\bar y+\bar z^3+\bar u^4=0\}/\ZZ_a(1,-b,-c,-d),\\
U_2=\{\bar x+\bar z^3+\bar u^4=0\}/\ZZ_b(-a,1,-c,-d),\\
U_3=\{\bar x\bar y+\bar u^{4}+1=0\}\ZZ_c(-a,-b,1,-d),\\
U_4=\{\bar x\bar y+\bar z^{3}+1=0\}\ZZ_d(-a,-b,-c,1).
\end{gather*}
Classification of terminal singularities (see \cite{YPG}) tells us that singularities
of $Y$ in charts $U_1$ and $U_2$ are terminal if one of the cases are realized
\begin{enumerate}
\item $c+d\equiv 0\mod a$ and $c+d\equiv 0\mod b$. This case is impossible since there
are $k,t\in \ZZ_+$ such that $c+d=ak$, $c+d=bt$ which bring us to the contradiction:
$b=\frac {12}7\left(k-\frac 7{12}\right)a$, $k=\frac{7t}{12t-7}<0$.

\item $c+d\equiv 0\mod a$ and $c\equiv 1\mod b$. This case is also impossible since there are
$k,t\in \ZZ_+$ such that $c+d=ak$, $c=bt+1$.We have either $t=1$ which brings us to the contradiction
with the condition $a+b=3c=4d$, or $t>1$. Therefore, we have $3tb+3=4ka-4tb-4$. Then, from the conditions
 $a=\frac {7tb+7}{4k}$ and $b(4k(3t-1)-7t)=7-12k$ we get  $b<0$ -- contradiction.

\item $c\equiv 1\mod a$ and $c\equiv 1\mod b$. In this case we have either $c=1$ which is impossible since
$a+b=3c=4d$ or $3c>a+b$ which produce the contradiction again.

\item $c\equiv 1\mod a$ and $d\equiv 1\mod b$. In this case we have got the contradiction with
$a+b=3c=4d$.
\end{enumerate}

\item Let's examine the terminality of $Y$. Indeed, $Y$ is covered  by four affine charts $U_1,\ U_2,\ U_3,\ U_4$
\begin{gather*}
U_1=\{\bar y+\bar z^3+\bar x^{4d-a-b}\bar u^4=0\}/\ZZ_a(1,-b,-c,-d),\\
U_2=\{\bar x+\bar z^3+\bar y^{4d-a-b}\bar u^4=0\}/\ZZ_b(-a,1,-c,-d),\\
U_3=\{\bar x\bar y+\bar z^{4d-a-b}\bar u^{4}+1=0\}\ZZ_c(-a,-b,1,-d),\\
U_4=\{\bar x\bar y+\bar z^{3}+\bar u^{4d-a-b}=0\}\ZZ_d(-a,-b,-c,1).
\end{gather*}
Classification of terminal singularities (see \cite{YPG}) tells us that singularities
of $Y$ in chart $U_4$ are terminal if $4d-a-b=1$ and singularities in charts $U_1$ and $U_2$ are terminal
if one fe the following cases are realized
\begin{enumerate}
\item $c+d\equiv 0\mod a$ and $c+d\equiv 0\mod b$. This case is impossible since
 there are $k,t\in \ZZ_+$ such that $c+d=ak$ and $c+d=bt$. We have $d=\frac {3ka+1}7$, $\frac {4a}t=b=\frac{4ka-1}{7t}$.
Then, we have got either $r<0$ or $a<0$ both lead us to the contradiction.

\item $c+d\equiv 0\mod a$ and $c\equiv 1\mod b$. This case is also impossible since there are
$k,t\in \ZZ_+$ such that $c+d=ak$, $c=bt+1$. We have either $t=1$ which contradicts to $a+b=3c=4d$ or
$t>1$. Therefore, $a=\frac {7bt+8}{4k}$ and $b(4k(3t-1)-7t)=8-12k$ - the contradiction with positivity of $b$.

\item $c+d\equiv 0\mod a$ and $d\equiv 1\mod b$. This case is impossible since there are $k,t\in \ZZ_+$
such that $c+d=ak$, $c=bt+1$. We have either $t=1$ -- the contradiction with
$a+b=3c=4d$ or $t>1$. Then, we have  $a=\frac {7bt+6}{3k}$ and $b(3k(4t-1)-7t)=6-9k$ -- the contradiction
with positivity of $b$.

\item $c\equiv 1\mod a$ and $c\equiv 1\mod b$. In this case we have got either  $c=1$ --
the contradiction with $a+b=3c=4d$ or $3c>a+b$ -- the contradiction again.

\item $c\equiv 1\mod a$ and $d\equiv 1\mod b$. This case leads us to the contradiction with
 $a+b=3c=4d$.
\end{enumerate}

\item The same calculations show us the impossibility of the last case.
\qed
\end{enumerate}

\subsection{Examination of the cases}
\subsubsection{Cases 2(a) and 2(b)}

Looking at the different models
we can reach the situation when the center of $S$ lies in the singularity of the model.
It follows from \cite{Kaw1}, that  $S$ can be realized
by some weighted blow-up $\bar Y\to X$. It is easy to notice that in this case $\bar Y \simeq Y$.
Lemma 3.2 leads us to the contradiction with a terminality of $Y$.

\subsubsection{Cases 2(c) and 2(d)}

If in the case 2(c) the center of $S$ lies in the singularity $A_2$ on $E_1$ then
$S$ can be realized by a weighted blow-up of $X$. It follows from the toric structure of weighted blow-ups and
from \cite{Kawk0} since the singularity (on $E_1$) $A_2$ lies in the origine of the chart $U_4$.
If the center of $S$ lies in another point then we can proceed in the same way as it was done
in the previous section.
The exact consequence
\begin{equation*}
0\longrightarrow\OO_{P^1}(-1)\longrightarrow\N_{X|\bar l}
\longrightarrow\OO_{P^1}(-1)\longrightarrow 0
\end{equation*}
allow us to blow-up  $\bar  l$ and then contract the surface to another ruling (see \cite{Kul}).
Therefore, we can consider the center of $S$ on some model to be in the  $x^2+y^2+z^2+u^2=0$.
It follows from \cite{Corti} that in this case $S$ realizes by the ordinary blow-up of this point.
Hence, in our model $S$ realizes as a blow-up of $\bar l$. Thus, the case 2(c) is not possible.
\par

The case 2(d) we will consider on the weighted blow-up of $X$ with weights $(1,2,1,1)$,
$E_1$ is the exceptional divisor.
It follows from \cite{Kaw1} that a divisorial contraction to a curve passing through
a terminal cyclic quotient singularity is a weighted blow-up of this singularity. Therefore, we get a
contradiction with lemma 2.5.
\qed

The main theorem of this section is proved


\begin{thebibliography}{99}
\bibitem{Corti}
\emph{Corti A.}
Singularities of linear systems and 3-fold
birational geometry, Explicit Birational Geometry of 3-folds,
A. Corti, M. Reid editors, L.M.S. lecture Note Series {\bf 281} (2000) 259--312

\bibitem{CorMe}
\emph{Corti A., Mella M.}
Birational geometry of terminal quartic 3-folds. I,
preprint, math AG/0102096, 2001

\bibitem{Cut}
\textit{Cutkosky S.} Elementary contractions of Gorenstein
threefolds, Math. Ann., \textbf{280} (1988), 521-525

\bibitem{Dan82}
\textit{Danilov V.I.} Birational geometry of toric 3-folds, Math.USSR Izv., \textbf{21} (1983),269-280

\bibitem{Fed2}
\textit{Fedorov I.Yu.} Blow-ups of three-dimensional terminal singularities: cA case,
math.AG/0111307

\bibitem{Kaw1}
\textit{Kawamata Y.} Divisorial contractions to 3-dimentional
terminal quotient singularities, Higher-dimentional complex
varieties (Trento, 1994), de Gruyter, 1996, 241-246

\bibitem{KMM}
\emph{Kawamata  Y.,Matsuda  K.,Matsuki  K.} Introduction  to  the
minimal model program //  Algebraic Geometry,Sendai. Adv. Stud. Pure Math. 1987.
{\bf 10.}  283--360

\bibitem{Kawk0}
\textit{Kawakita M.} Divisorial contractions in dimention three
which contract divisor to smooth points, preprint, math.
AG/0005207

\bibitem{Kawk1}
\textit{Kawakita M.} Divisorial contractions in dimension three
which contract divisors to compound $A_1$ points, preprint, math.
AG/0010207

\bibitem{Kawk2}
\textit{Kawakita M.} General elements in anticanonical systems of
threefolds with divisorial contractions and applications to classification, preprint,. math.AG/0110050

\bibitem{Kul}
\emph{Kulikov Vik. S.}
Degenerations of K3 surfaces and Enriques surfaces,Math. USSR Izv. \textbf{11} (1977), 957-989


\bibitem{Luo}
\textit{Luo T.} Divisorial contractions of threefolds:divisor to
point, Amer. J. Math., \textbf{120} (1998), 441-451


\bibitem{Mori85}
\textit{Mori S.} On 3-dimentional terminal singularities, Nagoya
Math. J.,\textbf{98} (1985), 43-66



\bibitem{Mori}
\textit{Mori S.} Threefolds whose canonical bundles are not
numerically effective, Ann. of Math., \textbf{116} (1982), 133-176


\bibitem{Ms84}
\textit{Morrison D., Stevens G.} Terminal quotient singularities
in dimention three and four, Proc. Amer. Math. Soc., \textbf{90}
(1984), 15-20


\bibitem{YPG}
\textit{Reid M.} Young person's guide to canonical singularities,
Proc. Symp. in Pure Math. \textbf{46} (1987) 343--416

\bibitem{Reid83}
\textit{Reid M.} Minimal models of canonical threefolds, Adv.
Stud. Pure Math. \textbf{1} (1983) 131--180


\bibitem{Sh1}
\emph{Shokurov V.V.} 3-folds log flips,Russian Acad. Sci. Izv. Math. \textbf{40},(1993), 95-202.

\end{thebibliography}
\end{document}